\documentclass[10pt]{amsart}
\usepackage{amsfonts}
\usepackage{ifthen}
\usepackage{amsthm}
\usepackage{amsmath}
\usepackage{graphicx}
\usepackage{amscd,amssymb,amsthm}

\newcounter{minutes}
\divide\time by 60
\newcounter{hours}
\multiply\time by 60 \addtocounter{minutes}{-\time}

\newcommand{\real}{\operatorname{Re}}

\newtheorem{theorem}{Theorem}
\newtheorem{corollary}{Corollary}
\newtheorem{lemma}{Lemma}
\newtheorem{definition}{Definition}


\title[Subordinations and superordinations for Bessel functions]{Differential subordinations and superordinations for generalized Bessel functions}

\author{Huda A. Al-Kharsani}
\address{Department of Mathematics, Girls College, University of Dammam, 838 Dammam, Saudi Arabia}
\email{halkharsani@ud.edu.sa}

\author{\'Arp\'ad Baricz}
\address{Department of Economics, Babe\c{s}-Bolyai University, 400591 Cluj-Napoca, Romania} \email{bariczocsi@yahoo.com}

\author{K.S. Nisar}
\address{Department of Mathematics, College of Arts and Science, Salman bin Abdulaziz University, 54 Wadi Addawasir, Saudi Arabia}
\email{ksnisar1@gmail.com}

\begin{document}

\def\thefootnote{}
\footnotetext{ \texttt{File:~\jobname .tex,
          printed: \number\year-0\number\month-\number\day,
          \thehours.\ifnum\theminutes<10{0}\fi\theminutes}
} \makeatletter\def\thefootnote{\@arabic\c@footnote}\makeatother

\maketitle

\begin{abstract}
Differential subordination and superordination preserving properties for
univalent functions in the open unit disk with an operator involving generalized
Bessel functions are derived. Some particular cases involving trigonometric
functions of our main results are also pointed out.
\end{abstract}

\section{\bf Introduction and some preliminary results}
\setcounter{equation}{0}

It is known that the generalized hypergeometric functions play
an important role in geometric function theory, especially in the
solution by de Branges of the famous Bieberbach conjecture. Motivated
by this, geometric properties (like univalence, starlikeness, convexity)
of different types of hypergeometric functions were investigated by many
authors. For example, Miller and Mocanu \cite{millerhyper} employed the
method of differential subordinations \cite{millerbook} to investigate the local univalence,
starlikeness and convexity of certain hypergeometric functions. Moreover, it is worth mentioning
that, motivated by the results of Miller and Mocanu, further results on
hypergeometric functions were obtained by Ponnusamy and Vuorinen \cite{pv1,pv2}.
Motivated by the above mentioned results some similar developments were also made for the
so-called generalized Bessel functions \cite{andras,bariczbook,bariczDebrecen,bariczthesis,bariczHardytr,bariczpuma,erhan,bariczponnusamy,mondal}. In this paper we make a
further contribution to the subject by showing some differential subordination and superordination
results for an operator involving the generalized Bessel functions of the first kind.

The generalized Bessel function of the first kind $w=w_{p,b,c}$ is
defined as the particular solution of the second-order linear
homogeneous differential equation \cite{bariczbook,bariczthesis}
\begin{equation}\label{eq:1.1} z^2
w''(z)+bzw'(z)+\left(cz^2-p^2+(1-b)p\right)w(z)=0,
\end{equation} which is natural
generalization of Bessel differential equation. This function has the representation
\begin{equation} \label{eq:1.1}
w(z)=w_{p,b,c}(z)=\sum_{n\geq0} \frac{(-1)^n c^n}{n! \Gamma
\left(p+n+\frac{b+1}{2}\right)}\left(\frac{z}{2}\right)^{2n+p},
\end{equation}
where $b,p,c,z\in\mathbb{C}$ and $c\neq0.$ The differential equation (1.1) permits the study of
Bessel, modified Bessel and spherical Bessel functions all together.
Solutions of (1.1) are referred to as the generalized Bessel
function of order $p.$ The particular solution given by (1.2) is
called the generalized Bessel function of the first kind of order
$p$. Although the series defined above is convergent everywhere, the
function $w_{p,b,c}$ is generally not univalent in the open unit disk $\mathbb{D}=\left\{z\in\mathbb{C}:|z|<1\right\}$. It
is worth mentioning that, in particular, when $b=c=1,$ we reobtain
the Bessel function of the first kind $w_{p,1,1}=J_{p}, \, \mbox { and for } \, c=-1
\,\mbox { and } \, b=1,$ the function $w_{p,1,-1}$ becomes the
modified Bessel function of the first kind $I_p$. Now, consider the function
$u_{p,b,c}:\mathbb{C}\to\mathbb{C},$ defined by the transformation
\[u_{p,b,c}(z)=2^{p}\Gamma \left (p+\frac{b+1}{2}\right) \cdot z^{-p/2} w_{p,b,c} (\sqrt{z}).\]
By using the well-known Pochhammer (or Appell ) symbol, defined in
terms of the Euler gamma function,
\[(a)_n=\frac{\Gamma (a+n)}{\Gamma (a)}=a(a+1)\ldots (a+n-1)\]
and $(a)_{0}=1,$ we obtain for the function $u_{p,b,c}$ the
following representation
\[\displaystyle u_{p,b,c}(z)=\sum_{n\geq0} \frac{(-c/4)^n}{{\left(p+\frac{b+1}{2}\right)}_{n}}\frac{z^n}{n!},\]
where $p+\frac{b+1}{2} \neq 0,-1,-2,{\dots}. $   This function is
analytic on $\mathbb{C}$ and satisfies the second order linear
differential equation
\[4z^2 u'' (z)+2(2p+b+1)zu' (z)+czu(z)=0.\]

Now, let $\mathcal{H}=\mathcal{H} (\mathbb{D})$ denote the class of
analytic functions defined in $\mathbb{D}$, and for $n \in
\mathbb{N}$ and $a \in \mathbb{C}$ let
$\mathcal{H}[a,n]$ consist of functions $f \in \mathcal{H}$ of the
form $f(z)=a+a_n z^n+a_{n+1}z^{n+1}+{\dots}.$ Let $ f $
and $F$ be members of $\mathcal{H}$. The function $f$ is said to be
subordinate to $F,$ or $F$ is said to be superordinate to $f$, if
there exists a function $w$ analytic in $\mathbb{D},$ with
$|w(z)|\leq |z|$ such that $f(z)=F(w(z)).$ In such a case, we
write $f\prec F $ or $f(z) \prec F(z).$ If the function $F$ is
univalent in $\mathbb{D}$, then $f\prec F$ if and only if $f(0)=F(0)
$ and $f(\mathbb{D})\subset F(\mathbb{D}) $ (cf. \cite{millerbook}). Let
$\varphi :\mathbb{C}^2\rightarrow \mathbb{C},$ and let $h$ be
univalent in $\mathbb{D}$. The subordination $\varphi
(p(z),zp'(z))\prec h(z)$ is called a first-order differential
subordination. Finally, let $\mathcal{A}$ denote the class of functions
\[f(z)=z+\sum_{n\geq1}a_{n+1}z^{n+1}\]\\
that are analytic and univalent in the open unit disk $\mathbb{D},$ and consider the following operator
$$B_{\kappa,c}(f)(z)=zu_{p,b,c}(z)* f(z)=z+\sum_{n\geq1}a_{n+1}
\frac{(-c/4)^n}{(\kappa)_{n}}\frac{z^{n+1}}{n!}$$
where $\kappa=p+\frac{b+1}{2}\neq 0,-1,-2,{\dots }.$ We mention that for this operator the next identity is valid
\begin{equation}
\label{eq:1.4}
z(B_{\kappa+2,c}(f(z)))'=(\kappa+1)B_{\kappa+1,c}(f(z))-\kappa B_{\kappa+2,c}(f(z)).
\end{equation}

In this paper some subordination and superordination preserving properties for
univalent functions in the open unit disk associated with the above operator will be derived.
The paper is organized as follows: this section contains the definitions and some preliminary results which will be used in the sequel. Section 2 contains the main results together with their consequences, while section 3 is devoted for the proofs of the main results. We note that very recently some other differential subordination and superordination results on the above operator were obtained in \cite{erhan}.

The following definitions and lemmas will be used in our
present investigation. For more details see \cite{millercomplvar,millerbook,millerbriot}.

\begin{definition} \cite{millerbook}
Let $\varphi:\mathbb{C}^2\rightarrow\mathbb{C},$ and let $h$ be univalent in
$\mathbb{D}$. If $p$ is analytic in $\mathbb{D}$ and satisfies the
differential subordination
\begin{equation} \label{eq:1.5} \varphi (p(z),zp'(z))\prec h(z),\end{equation}
then $p$ is called a solution of differential subordination \eqref{eq:1.5}. A
univalent function $q$ is called a dominant of the solutions of
differential subordination \eqref{eq:1.5}, or more simply a dominant, if
$p\prec q$ for all $p$ satisfying \eqref{eq:1.5}. A dominant $\tilde{q}$ that
satisfies $\tilde{q}\prec q$ for all dominants $q$ of \eqref{eq:1.5} is said
to be the best dominant of \eqref{eq:1.5}.
\end{definition}

\begin{definition} \cite{millercomplvar}
Let $\varphi:\mathbb{C}^2\rightarrow\mathbb{C},$ and let $h$ be analytic in
$\mathbb{D}.$ If $p$ and $z\mapsto \varphi(p(z),zp'(z))$ are univalent in
$\mathbb{D}$ and satisfy the differential superordination
\begin{equation}\label{eq:1.1} h(z)\prec \varphi (p(z), zp'(z)),\end{equation}
then $p$ is called a solution of differential superordination \eqref{eq:1.1}. An analytic
function $q$ is called a subordinant of the solutions of differential superordination \eqref{eq:1.1}, or more simply a subordinant, if
$q\prec p$ for all $p$ satisfying \eqref{eq:1.1}. A univalent subordinant
$\tilde{q}$ that satisfies $q\prec \tilde{q}$ for all subordinants
$q$ of \eqref{eq:1.1} is said to be the best subordinant \eqref{eq:1.1}.
\end{definition}

\begin{definition} \cite{millerbook}
Denote by $Q$ the class of functions $f$
that are analytic and injective on $\overline{\mathbb{D}}\backslash
E(f),$ where
\[E(f)=\left \{\zeta \in \partial \mathbb{D}:\lim_{z\rightarrow \zeta}f(z)=\infty \right\},\]
and $f$ is such that $f'(\zeta)\neq 0 \mbox{ for } \zeta \in \partial
\mathbb{D}\backslash E(f).$
\end{definition}

The following preliminary results which will be used in the sequel are some known results on admissible functions.

\begin{lemma} \cite{millerbook}\label{lem:1}
Suppose that the function
$\eta:\mathbb{C}^2\rightarrow \mathbb{C}$ satisfies the condition
\[\real\eta(\mathrm{i}s,t) \leq 0,\]
for all real $s$ and $t\leq -m(1+s^2)/2,$ where $m$ is a
positive integer. If the function $p(z)=1+p_{m}z^{m}+{\dots}$ is
analytic in $\mathbb{D}$ and
\[\real \eta(p(z),zp'(z))> 0\ \ \ \mbox{for all}\ \ \ z \in \mathbb{D}\]
then $\real p(z)> 0$ for all $z\in\mathbb{D}.$
\end{lemma}

\begin{lemma} \cite{ponnusamyroning} \label{lem:2}
Let $\beta,\gamma \in \mathbb{C}$ with $\beta\neq 0,$ and let $h\in \mathcal{H}(\mathbb{D})$ with
$h(0)=c.$ If $$\real(\beta h(z)+\gamma)> 0$$ for all $z\in\mathbb{D},$
then the solution of the differential equation
\[q(z)+\dfrac{zq'(z)}{\beta q (z)+\gamma}=h(z)\]
with $q(0)=c$ is analytic in $\mathbb{D}$ and satisfies $$\real(\beta q
(z)+\gamma)>0$$ for all $z\in\mathbb{D}.$
\end{lemma}

\begin{lemma}\cite{millerbook}\label{lem:3}
Let $p \in Q $ with $p(0)=a,$ and let
$q(z)=a+a_{n}z^{n}+{\dots}$ be analytic in $\mathbb{D}$ with
$q(z)\not\equiv a$ and $n\geq 1.$ If $q$ is not subordinate to $p$,
then there exist $z_0=r_0 e^{\mathrm{i}\theta} \in \mathbb{D}$ and
$\zeta_{0}\in \partial \mathbb{D}\backslash E(p), $ for which
$q(r_{0}\mathbb{D})\subset p(\mathbb{D}),$ and
\[q(z_0)=p(\zeta_{0}), \quad \quad
z_{0}q'(z_0)=m\zeta_{0}p'(\zeta_{0}) \quad (m\geq n).\]
\end{lemma}

We will use also the concept Loewner subordination chain, which is defined as follows.

\begin{definition} \cite{millercomplvar}
A function $(z,t)\mapsto L(z,t)$ defined on $\mathbb{D}\times [0,\infty)$ is a
subordination chain (or Loewner chain) if $L(\cdot, t)$ is analytic
and univalent in $\mathbb{D}$ for all $t \in [0,\infty)$ and $L(z,\cdot)$
is continuously differentiable on $[0,\infty)$ for all
$z\in\mathbb{D},$ and $L(z,s)\prec L(z,t) $ for $0\leq s< t.$
\end{definition}

The next results are also useful in order to obtain the main results of this paper.

\begin{lemma} \cite{millercomplvar}\label{lem:4}  Let $q \in \mathcal{H}[a,1],$ $\varphi
:\mathbb{C}^2\rightarrow \mathbb{C},$ and set $h(z)\equiv
\varphi(q(z),zq'(z)).$ If $$L(z,t)=\varphi (q(z),tzq'(z))$$ is a
subordination chain and $p \in \mathcal{H}[a,1] \cap Q,$ then
\[h(z)\prec \varphi (p(z),zp'(z))\]
implies that
\[q(z)\prec p(z).\]
Furthermore, if $\varphi (q(z),zq'(z))=h(z)$ has a univalent
solution $q\in Q$, then $q$ is the best subordinant.
\end{lemma}

\begin{lemma} \cite{millercomplvar}\label{lem:5}
The function $(z,t)\mapsto L(z,t)=a_1(t)z+{\dots},$ with
$a_1(t)\neq 0$ and $\displaystyle \lim_{t\rightarrow
\infty}|a_{1}(t)|=\infty,$ is
a subordination chain if and only if
\[\real \left(\frac{z\partial L (z,t)/ \partial z}{\partial L(z,t)/\partial t}\right)> 0\ \ \  \mbox{for all}\ \ \ z \in \mathbb{D},\ 0 \leq t <\infty.\]
\end{lemma}

\section{\bf Main results and their consequences}
\setcounter{equation}{0}

Our first main result is the following theorem.

\begin{theorem}\label{th1}
Let $f,g \in \mathcal{A},$ $\lambda\in[0,1),$ $c\in\mathbb{C}$ with $c\neq0$ and $p,b\in\mathbb{R}$ be such that $\kappa>-1.$ Let also
$$\Phi(z)=(1-\lambda) \frac{B_{\kappa+1,c}(g(z))}{z} + \lambda \frac{B_{\kappa+2,c}(g(z))}{z},$$
and suppose that
\begin{equation} \label{eq:2.1} \real \left(\frac{z\Phi''(z)}{\Phi'(z)}+1\right)>
-\gamma_{\lambda,\kappa} \ \ \ \mbox{for all}\ \  z \in \mathbb{D},\end{equation}
where
$$\gamma_{\lambda,\kappa}=\frac{(1-\lambda)^2+(\kappa+1)^2-\sqrt{(1-\lambda)^4+(\kappa+1)^4}}{4(1-\lambda)(\kappa+1)}.$$
Then the subordination condition
$$(1-\lambda)\frac{B_{\kappa+1,c}(f(z))}{z}+\lambda
\frac{B_{\kappa+2,c}(f(z))}{z} \prec (1-\lambda) \frac{B_{\kappa+1,c}(g(z))}{z} + \lambda \frac{B_{\kappa+2,c}(g(z))}{z}$$
implies that
$$\frac{B_{\kappa+2,c}(f(z))}{z}\prec \frac{B_{\kappa+2,c}(g(z))}{z}.$$
Moreover, the function $z\mapsto B_{\kappa+2,c}(g(z))/z$ is the best dominant.
\end{theorem}

Now, choosing $\lambda=0$ in the above theorem, we have the following result.

\begin{corollary}
Let $c\in\mathbb{C}$ with $c\neq0,$ $f,g\in\mathcal{A}$ and suppose that $p,b\in\mathbb{R}$ are such that $\kappa>-1.$ Consider also the function $\Psi:\mathbb{D}\to\mathbb{C},$ defined by $\Psi(z)={B_{\kappa+1,c}(g)(z)}/{z},$ and suppose that the condition $$\real\left(\frac{z\Psi''(z)}{\Psi'(z)}+1\right)>\frac{1+(\kappa+1)^2-\sqrt{1+(\kappa+1)^4}}{4(\kappa+1)}$$ is satisfied for all $z\in\mathbb{D}.$ Then
$$\frac{B_{k+1,c}(f)(z)}{z}\prec\frac{B_{k+1,c}(g)(z)}{z}$$ implies
that
$$\frac{B_{k+2,c}(f)(z)}{z}\prec\frac{B_{k+2,c}(g)(z)}{z}.$$
Moreover the function $z\mapsto{B_{k+2,c}(g)(z)}/{z}$ is the best dominant.
\end{corollary}

Taking into account the above results, we have the following
particular cases. Choosing $f(z)=\dfrac{z}{1-z}$ and $g(z)=z+az^{2}$
where $|a|<\frac{1}{2}$ in the above corollary we obtain that for $\kappa>-1$ the subordination
$$u_{p+1,b,c}(z) \prec 1-\frac{acz}{4(\kappa+1)}$$
implies that $$u_{p+2,b,c}(z) \prec 1-\frac{acz}{4(\kappa+2)}$$
or equivalently
$$ \left|u_{p+1,b,c}(z)-1\right|<\left|\frac{ac}{4(\kappa+1)}\right|\ \ \Rightarrow \left|u_{p+2,b,c}(z)-1\right|<
\left|\frac{ac}{4(\kappa+2)}\right|.$$
It is important to note here that the above result is related to a recent open problem from \cite{andras} concerning a subordination property of normalized Bessel functions with different parameters. Now, choosing in the above inequalities $p=-\frac{3}{2},$ $b=c=1$ ($\kappa=-\frac{1}{2}$) and $p=-\frac{1}{2},$ $b=c=1$ ($\kappa=\frac{1}{2}$), respectively, we obtain for all $z\in\mathbb{D}$ and $|a|<\frac{1}{2}$ the following chain of implications
$$\left|\cos\sqrt{z}-1\right|< \frac{|a|}{2} \Rightarrow \left|\frac{\sin\sqrt{z}}{\sqrt{z}}-1\right|<\frac{|a|}{6} \Rightarrow \left|3\frac{\sin\sqrt{z}}{z\sqrt{z}}-3\frac{\cos\sqrt{z}}{z}-1\right|<\frac{|a|}{10}.$$ Here we used the relations
$$u_{-\frac{1}{2},1,1}(z)=\sqrt{\frac{\pi}{2}}z^{\frac{1}{4}}J_{-\frac{1}{2}}(\sqrt{z})=\cos\sqrt{z},$$
$$u_{\frac{1}{2},1,1}(z)=\sqrt{\frac{\pi}{2}}z^{-\frac{1}{4}}J_{\frac{1}{2}}(\sqrt{z})=\frac{\sin\sqrt{z}}{\sqrt{z}}$$
and
$$u_{\frac{3}{2},1,1}(z)=3\sqrt{\frac{\pi}{2}}z^{-\frac{3}{4}}J_{\frac{3}{2}}(\sqrt{z})=3\left(\frac{\sin\sqrt{z}}{z\sqrt{z}}-\frac{\cos\sqrt{z}}{z}\right).$$

Now we consider the dual of Theorem \ref{th1}.

\begin{theorem}\label{th2}
Let $f,g \in \mathcal{A},$ $p,b\in\mathbb{R}$ such that $\kappa>-1,$ $c\in\mathbb{C}$ with $c\neq0$ and $\lambda\in[0,1).$ Let also
\[\Phi(z)=(1-\lambda)\frac{B_{\kappa+1,c}(g(z))}{z}+\lambda \frac{B_{\kappa+2,c}(g(z))}{z}.\]
Suppose that for all $z\in\mathbb{D}$ we have
\[\real \left(\frac{z\Phi''(z)}{\Phi'(z)}+1\right)>-\gamma_{\lambda,\kappa}\]
and assume that
$$z\mapsto (1-\lambda)\frac{B_{\kappa+1,c}(f(z))}{z}+\lambda
\frac{B_{\kappa+1,c}(f(z))}{z}$$ is univalent in $\mathbb{D},$ and $z\mapsto B_{\kappa+2}(f(z))/
z \in \mathcal{H}[1,1]\cap Q.$ Then the superordination
\begin{equation} \label{eq:2.5ii}
(1-\lambda)\frac{B_{\kappa+1,c}(g(z))}{z}+\lambda \frac{B_{\kappa+2,c}(g(z))}{z}\prec
(1-\lambda)\frac{B_{\kappa+1,c}(f(z))}{z}+\lambda \frac{B_{\kappa+2,}(f(z))}{z}\end{equation} implies that
$$\frac{B_{\kappa+2}(g(z))}{z}\prec \frac{B_{\kappa+2}(f(z))}{z}.$$ Moreover, the function $z\mapsto B_{\kappa+2}(g(z))/z$ is the best
dominant.
\end{theorem}

Combining Theorems \ref{th1} and \ref{th2} we get the following sandwich type result.

\begin{corollary}\label{th3}
Let $f,g_1,g_2\in \mathcal{A}$ and suppose that $p,b\in\mathbb{R}$ such that $\kappa>-1,$ $\lambda\in[0,1),$ $c\in\mathbb{C}$ with $c\neq0.$ Consider the functions
$\Phi_1,\Phi_2:\mathbb{D}\to\mathbb{C},$ defined by
$$\Phi_{i}(z)=(1-\lambda)\frac{B_{\kappa+1,c}(g_{i}\,(z))}{z}+\lambda
\frac{B_{\kappa+2,c}(g_{i}\,(z))}{z}, \ \ \ i=1,2,$$ and suppose that for $i=1,2$ and $z\in\mathbb{D}$ we have
$$\real\left(1+\frac{z\Phi''_{i}(z)}{\Phi'_i(z)}\right)>
-\gamma_{\lambda,\kappa}.$$ Moreover, assume that
\[z\mapsto (1-\lambda)\frac{B_{\kappa+1,c}(f(z))}{z}+\lambda \frac{B_{\kappa+2,c}(f(z))}{z}\]
is univalent in $\mathbb{D}$ and $z\mapsto B_{\kappa+2,c}(f(z))/ z \in \mathcal{H}[1,1]\cap
Q.$ Then
\[\Phi_{1}(z)\prec (1-\lambda)\frac{B_{\kappa+1,c}(f(z))}{z}+\lambda \frac{B_{\kappa+2,c}(f(z))}{z}\prec \Phi_{2}(z)\]
implies that
\[\frac{B_{\kappa+2,c}(g_{1}(z))}{z}\prec \frac{B_{\kappa+2,c}(f(z))}{z}\prec\frac{B_{\kappa+2,c}(g_{2}(z))}{z}.\]
Moreover, the function $z\mapsto {B_{\kappa+2,c}(g(z))}/{z}$ and
$z\mapsto {B_{\kappa+2,c}(g_{2}(z))}/{z}$ are the best subordination and the
best dominant, respectively.
\end{corollary}

Finally, let us consider the generalized Libera integral operator
$$F_{\mu}(f)(z)=\frac{\mu +1}{z^\mu} \int^{z}_{0}
t^{\mu -1} f(t)dt,$$
where $\mu > -1$ and $f\in\mathcal{A}.$ The following theorem is a sandwich-type result involving the
generalized Libera integral operator $F_{\mu}(f).$

\begin{theorem}\label{libera}
Let $b,p,c\in\mathbb{C}$ with $c\neq0,$ $f,g_1,g_2\in \mathcal{A}$ and let $\omega_{i}(z)={B_{\kappa,c}(g_i(z))}/{z}$ for $i=1,2.$
Suppose that $\mu>-1$ and for $z\in\mathbb{D}$ we have
$$\real \left(\frac{z\omega''_{i}(z)}{\omega_{i}'(z)}+1\right)>
-\gamma_{\mu},$$ where
$$\gamma_{\mu}=\frac{1+(\mu+1)^2-\sqrt{1+(\mu+1)^4}}{4(\mu+1)}.$$
If $z\mapsto {B_{\kappa,c}(f(z))}/{z}$ is univalent in $\mathbb{D}$ and
$z\mapsto B_{\kappa,c}(F_{\mu}(f))(z) \in \mathcal{H}[1,1]\cap Q, $ then
\[\omega_{1}(z)\prec \frac{B_{\kappa,c}(f(z))}{z}\prec \omega_{2}(z) \]
implies that
\[\frac{B_{\kappa,c}(F_{\mu}(g_1))(z)}{z}\prec \frac{B_{\kappa,c}(F_{\mu}(f))(z)}{z}\prec
\frac{B_{\kappa,c}(F_{\mu}(g_2))(z)}{z}.\]
Moreover, the functions $z\mapsto {B_{\kappa,c}(F_{\mu}(g_1))(z)}/{z}$ and $z\mapsto {B_{\kappa,c}(F_{\mu}(g_2))(z)}/{z}$ are the best
subordinant and the best dominant, respectively.
\end{theorem}

\section{\bf Proofs of the Main Results}
\setcounter{equation}{0}

In this section our aim is to present the proofs of the main results.

\begin{proof}[\bf Proof of Theorem \ref{th1}]
Let us define the functions $\phi,\psi:\mathbb{D}\to\mathbb{C}$ by
\begin{equation} \label{eq:2.4} \phi (z) =\frac{B_{\kappa+2,c}(g(z))}{z}, \quad
\quad \psi(z)=\frac{B_{\kappa+2,c}(f(z))}{z}.\end{equation}
We first show that \eqref{eq:2.1} implies that for all $z\in\mathbb{D}$ we have
\begin{equation} \label{eq:2.1i} \real\left(\frac{z\phi''(z)}{\phi'(z)}+1\right)> 0.\end{equation}
Differentiating both sides of the first equation in \eqref{eq:2.4} and using \eqref{eq:1.4} for $g\in\mathcal{A}$ we obtain
\begin{equation} \label{eq:2.4i}
\frac{\kappa+1}{1-\lambda}\Phi(z)=z\phi'(z)+\frac{\kappa+1}{1-\lambda}\phi(z).\end{equation}
Now, differentiating twice both sides of \eqref{eq:2.4i} yields the following
$$1+\frac{z\Phi''(z)}{\Phi'(z)}=q(z)+\frac{zq'(z)}{q(z)+\frac{\kappa+1}{1-\lambda}}\equiv
h(z),\ \ \mbox{where}\ \ q(z)=1+\frac{z\phi''(z)}{\phi'(z)}.$$
We note that from \eqref{eq:2.1} we obtain that for all $z\in\mathbb{D},$ $\lambda\in[0,1)$ and $\kappa>-1$
\[\real \left(h(z)+\frac{\kappa+1}{1-\lambda}\right)>\real \left(h(z)+\gamma_{\lambda,\kappa}\right) > 0,\]
where we used the inequality
$$-\frac{(\kappa+1)^4}{(1-\lambda)^2+\sqrt{(1-\lambda)^4+(\kappa+1)^4}}=(1-\lambda)^4-\sqrt{(1-\lambda)^4+(\kappa+1)^4}<3(\kappa+1)^2.$$
By using Lemma \ref{lem:2}, we conclude that the differential equation
$$q(z)+\frac{zq'(z)}{q(z)+\frac{\kappa+1}{1-\lambda}}=h(z)$$
has a solution $q \in \mathcal{H}(\mathbb{D})$ such that $q(0)=h(0)=1$. Now, let us consider the expression
$$\xi(u,v)=u+\frac{v}{u+\frac{\kappa+1}{1-\lambda}}+\gamma_{\lambda,\kappa}.$$
From \eqref{eq:2.1} it follows that for all $z\in\mathbb{D}$ we have $$\real \xi(q(z),zq'(z))> 0.$$
Now, to prove \eqref{eq:2.1i} we shall apply Lemma \ref{lem:1}. Thus, we need to show that $\real \xi(\mathrm{i}s,t) \leq 0$ for all real $s$
and $t\leq -\frac{1+s^2}{2}.$  But, this is true since for such $s$ and $t$ we have
\begin{align*}\real\xi(\mathrm{i}s,t)&=\dfrac{t\frac{\kappa+1}{1-\lambda}}{\left|\mathrm{i}s+\frac{\kappa+1}{1-\lambda}\right|^2}
+\frac{(1-\lambda)(\kappa+1)}{2\left((1-\lambda)^{2}+(\kappa+1)^{2}+\sqrt{(1-\lambda)^{4}+(\kappa+1)^{4}}\right)}\\
&\leq\frac{-\frac{s^{2}+1}{2}\frac{\kappa+1}{1-\lambda}}{s^{2}+\left(\frac{\kappa+1}{1-\lambda}\right)^{2}}+
\frac{\frac{\kappa+1}{1-\lambda}}{2\left(1+\left(\frac{\kappa+1}{1-\lambda}\right)^{2}+\sqrt{1+\left(\frac{\kappa+1}{1-\lambda}\right)^{4}}\right)}\\
&=-\frac{\kappa+1}{1-\lambda}\cdot\frac{1+s^{2}\left(\frac{\kappa+1}{1-\lambda}\right)^{2}+\sqrt{1+\left(\frac{\kappa+1}{1-\lambda}\right)^{4}}
+s^{2}\sqrt{1+\left(\frac{k+1}{1-\lambda}\right)^{4}}}{\left(s^{2}+\left(\frac{\kappa+1}{1-\lambda}\right)^{2}\right)
\left(1+\left(\frac{\kappa+1}{1-\lambda}\right)^{2}+\sqrt{1+\left(\frac{\kappa+1}{1-\lambda}\right)^{4}}\right)}\leq0.
\end{align*}
Applying Lemma \ref{lem:1} we obtain that $\real q(z)>0$ for all $z\in\mathbb{D},$ that is, indeed \eqref{eq:2.1i} is valid for all $z\in\mathbb{D}.$
In other, words the function $\phi $ is convex in $\mathbb{D}.$ Next we prove that
subordination condition of this theorem implies that $\psi\prec\phi.$ Without loss of generality, we can assume that $\phi$ is
analytic and univalent on $\overline{\mathbb{D}}$ and $\phi'(\xi)\neq 0$ for $|\xi|=1.$ For this purpose, we consider the function
$(z,t)\mapsto L(z,t)$ defined by
\[L(z,t)=\phi (z)+(1+t)\frac{1-\lambda}{\kappa+1}z\phi'(z),\]
where $z\in\mathbb{D}$ and $t\geq0.$ We note that for $\kappa>-1,$ $\lambda\in(0,1)$ and $t\geq0$ we have
\[\left.\frac{\partial L(z,t)}{\partial z}\right|_{z=0}=\phi'(0) \left(1+(1+t)\frac{1-\lambda}{\kappa+1}\right)\neq 0,\]
which shows that the function $(z,t)\mapsto L(z,t)=a_1(t)z+{\dots}$ satisfies the condition $a_1(t)\neq 0$ for all $t\in[0,\infty).$ Moreover,
for $\kappa>-1,$ $\lambda\in(0,1)$ and $t\geq0$ we have
$$\real \left(\dfrac{z\partial L(z,t)/ \partial z}{ \partial L (z,t)/ \partial t}\right)=\real\left(\frac{\kappa+1}{1-\lambda}+(1+t)\left(1+\frac{z\phi''(z)}{\phi'(z)}\right)\right)>0,$$
which by means of Lemma \ref{lem:5} shows that $(z,t)\mapsto L(z,t)$ is a subordination chain. We observe from the definition of a subordination chain that
$L(z,0)\prec L(z,t)$ for $t\geq0$ and hence $L(\varsigma,t)\not \in L(\mathbb{D},0)=\Phi(\mathbb{D})$ for $\varsigma \in \partial\mathbb{D}$ and $t\geq0.$ Now, suppose that $\psi$ is not subordinate to $\phi.$ Then by Lemma \ref{lem:3} there exist $z_{0}\in\mathbb{D}$ and $\xi_{0}\in\partial\mathbb{D}$ such that
\[\psi(z_{0})=\phi (\xi_{0}), \,\, z_{0}\psi'(z_{0})=(1+t)\xi_{0}\phi'(\xi_{0})\]
for $t\geq0.$ Hence \begin{align*}L(\xi_{0},t)&=\phi(\xi_{0})+(1+t)\frac{1-\lambda}{\kappa+1}\xi_{0}\phi'(\xi_{0})=\psi (z_{0})+\frac{1-\lambda}{\kappa+1}z_{0}\psi '(z_{0})\\&=(1-\lambda)\frac{B_{\kappa+1,c}(f(z_{0}))}{z_{0}}+\lambda\frac{B_{\kappa+2,c}(f(z_{0}))}{z_{0}}\in \Phi(\mathbb{D}),\end{align*}
which is a contradiction. Therefore, the subordination condition of the theorem must imply the subordination $\psi\prec\phi.$ Now, considering $\psi=\phi$ we can see that $\phi$ is the best dominant. This completes the proof of this theorem.
\end{proof}

\begin{proof}[\bf Proof of Theorem \ref{th2}]
The first part of the proof is similar to that of the proof of Theorem \ref{th1} and because of this
we will use the same notation as in the above proof. For this we define the functions $\phi$ and $\psi$ as in \eqref{eq:2.4} and consider the equation \eqref{eq:2.4i}, that is,
$$\Phi(z)=\phi(z)+\frac{1-\lambda}{\kappa+1}z\phi'(z).$$
This yields the relationship
$$1+\frac{z\Phi''(z)}{\Phi'(z)}=q(z)+\frac{zq'(z)}{q(z)+\frac{\kappa+1}{1-\lambda}}$$
and by using the same method as in the proof of Theorem \ref{th1}, we prove that $\real q(z)>0$ for all $z\in\mathbb{D}.$ That is, $\phi $ is convex and hence univalent in $\mathbb{D}$. Next, we prove that the superordination condition \eqref{eq:2.5ii} implies that $\phi\prec \psi.$ For this consider the function $(z,t)\mapsto L(z,t)$ defined by
$$L(z,t)=\phi (z)+\frac{1-\lambda}{\kappa+1}tz\phi'(z),$$
where $z\in\mathbb{D}$ and $0\leq t<\infty.$ Proceeding as in the proof
of Theorem \ref{th1} we can prove that $(z,t)\mapsto L(z,t)$ is a subordination chain.
Therefore according to Lemma \ref{lem:4} we conclude that superordination condition
\eqref{eq:2.5ii} must imply the superordination $\phi\prec \psi.$ Furthermore, since
the differential equation \eqref{eq:2.4i} has the univalent solution $\phi,$ it is the
best subordinant of the differential superordination.
\end{proof}

\begin{proof}[\bf Proof of Theorem \ref{libera}]
Let us define the functions $\psi,\phi_1,\phi_2:\mathbb{D}\to\mathbb{C}$ by
$$\psi(z)=\frac{B_{\kappa,c}(F_{\mu}(f))(z)}{z}, \quad
\quad \phi_{i}(z)=\frac{B_{\kappa,c}(F_{\mu}(g_i))(z)}{z}, \ i=1,2,$$
respectively. From the definition of the integral operator $F_{\mu}(f)$ we get
$$F_{\mu}(f)(z)=z+a_2\frac{\mu+1}{\mu+2}z^2+{\dots}+a_n\frac{\mu+1}{\mu+n}z^n+{\dots},$$
which implies that
$$z(B_{\kappa,c}(F_{\mu}(f))(z))'=(\mu+1)B_{\kappa,c}(f)(z)-\mu
B_{\kappa,c}(F_{\mu}(f))(z)=z+\sum_{n\geq 1}a_{n+1}\frac{\mu+1}{\mu+n+1}\frac{(-c/4)^{n}}{(\kappa)_n}\frac{z^{n+1}}{n!}.$$
Consequently, for $i=1,2$ we have
\[(\mu+1)\omega_{i}(z)=(\mu+1)\phi_{i}(z)+z\phi'_{i}(z)\]
and for $i=1,2$ using the notation
\[q_{i}(z)=1+\dfrac{z\phi''_{i}(z)}{\phi'_{i}(z)}\]
after differentiation we get for $i=1,2$
\[1+\frac{z\omega''_{i}(z)}{\omega'_{i}(z)}=q_{i}(z)+\dfrac{zq'_{i}(z)}{q_{i}(z)+\mu+1}\]
The remaining part of the proof is similar to that of the proof of Theorem \ref{th1}, and thus we omit the details.
\end{proof}

\subsection*{Acknowledgements} The research of H.A. Al-kharsani and K.S. Nisar was supported by Deanship of Scientific Research in University of Dammam, project ID$\#$2014219. The work of \'A. Baricz was supported by a research grant of the Romanian
National Authority for Scientific Research, CNCS-UEFISCDI, project number PN-II-RU-TE-2012-3-0190/2014.

\end{document}